\documentclass[12pt,twoside]{article}
\usepackage{amsmath, amssymb, amsthm,graphicx}
\font\fourteenb=cmb10 at 14pt
\begin{document}
\vspace*{-1.0cm}\noindent \copyright
 Journal of Technical University at Plovdiv\\[-0.0mm]\
\ Fundamental Sciences and Applications, Vol. 1, 1995\\[-0.0mm]
\textit{Series A-Pure and Applied Mathematics}\\[-0.0mm]
\ Bulgaria, ISSN 1310-8271\\[+1.2cm]
\font\fourteenb=cmb10 at 12pt
\begin{center}

{\bf \LARGE BMOA estimates and radial growth of  $ {B_{\phi} } $
functions
   \\ \ \\ \large Peyo Stoilov, Roumyana Gesheva, Milena Racheva }
\end{center}


\begin{abstract}
BMO estimates and the radial growth of Bloch functions have been
studied   by   B.   Korenblum [3]. The   present   paper contains
some   natural generalizations of these results.
\end{abstract}

     \section{Introduction}
\par
Let  $D$  denote the unit disk   $\left\{z\in
C:\left|z\right|<1\right\}$    and  $T$ - the unit circle
$\left\{z:\left|z\right|=1\right\}.$ The  space BMOA is the space
of functions  $f\in H^{1} $  for which
$$\left\| f\right\| _{*} =\mathop{\sup \frac{1}{m\left(I\right)}
\int \limits _{I}\left|f-f_{I} \right|dm }\limits_{} <\infty
,{\kern 1pt}{\kern 1pt} {\kern 1pt} {\kern 1pt} {\kern 1pt}where$$
$$f_{I} =\frac{1}{m\left(I\right)} \int \limits _{I}fdm
{\kern 1pt} ,{\kern 1pt} {\kern 1pt} {\kern 1pt} {\kern 1pt}
{\kern 1pt} {\kern 1pt} {\kern 1pt} {\kern 1pt} {\kern 1pt} {\kern
1pt} {\kern 1pt} {\kern 1pt} {\kern 1pt} {\kern 1pt} {\kern 1pt}
{\kern 1pt} I\subseteq T.$$
\par
Here  $m$  is normalized Lebesgue measure on   $T.$

     It is known that for an  analytic functions  $f$  in  $D$  the following conditions are equivalent (see, for example, [1] or [2]):

\

${\it a) }f\in BMOA ;{\kern 1pt}{\kern 1pt} {\kern 1pt} {\kern
1pt} {\kern 1pt}{\kern 1pt}{\kern 1pt} {\kern 1pt} {\kern 1pt}
{\kern 1pt}$
$${\it b)\left\| f\right\| _{BMOA}^{2} =\mathop{\sup
}\limits_{\xi \in D} \iint \limits
_{D}\left|f'\left(z\right)\right|^{2} \frac{(1-\left|z\right|^{2}
){\kern 1pt} {\kern 1pt} (1-\left|\xi \right|^{2} {\kern 1pt}
){\kern 1pt} {\kern 1pt} }{\left|1-\mathop{\xi }\limits^{\_ }
z\right|^{2} } dm_{2} <\infty ,{\kern 1pt}{\kern 1pt} {\kern 1pt}
{\kern 1pt} {\kern 1pt}{\kern 1pt}{\kern 1pt} {\kern 1pt} {\kern
1pt} {\kern 1pt}{\kern 1pt}{\kern 1pt} {\kern 1pt} {\kern 1pt}
{\kern 1pt}{\kern 1pt}{\kern 1pt} {\kern 1pt} {\kern 1pt} {\kern
1pt}{\kern 1pt}{\kern 1pt} {\kern 1pt} {\kern 1pt} {\kern
1pt}{\kern 1pt}{\kern 1pt} {\kern 1pt} {\kern 1pt} {\kern
1pt}{\kern 1pt}{\kern 1pt} {\kern 1pt} {\kern 1pt} {\kern
1pt}{\kern 1pt}{\kern 1pt} {\kern 1pt} {\kern 1pt} {\kern 1pt}}$$

where   $m_{2} $  denotes normalized  Lebesgue measure in  $D.$

     Let  $\phi \left(t\right)$  be a positive and continuous function on  $\left(0,1\right).$

     Let  $B_{\phi } $  denotes the space of all analytic functions  $f$  in  $D$  satisfying the condition

$$\left\| f\right\| _{B_{\phi } } =\mathop{\sup }\limits_{z\in D} \frac{(1-\left|z\right|^{2} )\left|f'(z)\right|}{\phi {\kern 1pt} {\kern 1pt} (1-\left|z\right|^{2} {\kern 1pt} )} <\infty .$$

\

     For     $\phi \left(t\right)=t^{\alpha } ,{\kern 1pt} {\kern 1pt} {\kern 1pt} {\kern 1pt} 0<\alpha \le 1,{\kern 1pt} {\kern 1pt} {\kern 1pt} {\kern 1pt} {\kern 1pt} {\kern 1pt} B_{\phi } =\Lambda _{\alpha } $   is the usual  Lipschitz  class.

     In the  case   $\alpha =0,$   $\Lambda _{0} $   is the Bloch space (usually denoted  $B$  ).

In this paper some results of   B. Korenblum  for  the Bloch space
$B$  are generalized  for  $B_{\phi }.$
\par
Note that  the  function $\log (1-z){\kern 1pt} {\kern 1pt} {\kern
1pt} {\kern 1pt} \in {\kern 1pt} {\kern 1pt} {\kern 1pt} B$ ,
however   $\log ^{2} (1-z){\kern 1pt} {\kern 1pt} {\kern 1pt}
\notin {\kern 1pt} {\kern 1pt} B.$
\section{BMOA estimates for  $B_{\phi } $  functions and {\kern 1pt} {\kern 1pt} {\kern 1pt}{\kern 1pt} {\kern 1pt} {\kern 1pt}applications}

\

     Let    $\phi (t)$  satisfies the condition

     $$\int \limits _{x}^{1}\frac{\phi ^{2} ({\kern 1pt} t)}{t} dt =g(x)<\infty {\kern 1pt} {\kern 1pt} {\kern 1pt} {\kern 1pt} {\kern 1pt} {\kern 1pt} for{\kern 1pt} {\kern 1pt} {\kern 1pt} all{\kern 1pt} {\kern 1pt} {\kern 1pt} {\kern 1pt} {\kern 1pt} 0<x<1.$$

       If  $f\in {\kern 1pt} {\kern 1pt} {\kern 1pt} B_{\phi } $   we write   $f_{r} (z){\kern 1pt} {\kern 1pt} {\kern 1pt} \mathop{=}\limits^{def} {\kern 1pt} {\kern 1pt} f(rz)$   for    $0<r<1$ .

\

     {\bf Theorem 1.} {\it Let  } $f\in {\kern 1pt} {\kern 1pt} {\kern 1pt} B_{\phi } $ {\it . Then}

     \

\begin{equation}
                                 \left\| f_{r} \right\|
_{BMOA} \le \left\| f\right\| _{B_{\phi } } \sqrt{g(1-r^{2} )}
{\kern 1pt} {\kern 1pt} {\kern 1pt} {\kern 1pt} ,{\kern 1pt}
{\kern 1pt} {\kern 1pt} {\kern 1pt} {\kern 1pt} {\kern 1pt} {\kern
1pt} {\kern 1pt} 0<r<1.
\end{equation} \

\

     {\it Proof.}  Let    $\xi \in D$ .  Then

\

     $${\kern 1pt} {\kern 1pt} \iint \limits _{D}\left|f'_{r} \left(z\right)\right|^{2} {\kern 1pt} {\kern 1pt} {\kern 1pt} \frac{(1-\left|z\right|^{2} ){\kern 1pt} {\kern 1pt} (1-\left|\xi \right|^{2} ){\kern 1pt} {\kern 1pt} {\kern 1pt} }{\left|1-\mathop{\xi }\limits^{\_ } z\right|^{2} } {\kern 1pt} {\kern 1pt} {\kern 1pt} dm_{2} = $$

\

     $$={\kern 1pt} {\kern 1pt} {\kern 1pt} {\kern 1pt} r^{2} \iint \limits _{D}\frac{\left|f'(rz)\right|^{2} (1-\left|rz\right|^{2} )^{2} }{\phi ^{2} (1-\left|rz\right|^{2} )} \cdot \frac{\phi ^{2} (1-\left|rz\right|^{2} )}{(1-\left|rz\right|^{2} )^{2} }  \cdot \frac{(1-\left|z\right|^{2} ){\kern 1pt} {\kern 1pt} (1-\left|\xi \right|^{2} {\kern 1pt} )}{\left|1-\mathop{\xi }\limits^{\_ } z\right|^{2} } dm_{2} \le $$

\

$$\le r^{2} \left\| f\right\| _{B_{\phi } }^{2} \iint \limits _{D}\frac{\phi ^{2} (1-\left|rz\right|^{2} )}{1-\left|rz\right|^{2} } \cdot \frac{1-\left|\xi \right|^{2} }{\left|1-\overline{\xi }{\kern 1pt} {\kern 1pt} z\right|^{2} } dm_{2} = $$

\

$$=2r^{2} \left\| f\right\| _{B_{\phi } }^{2} \int \limits _{0}^{1}\int \limits _{T}\frac{\phi ^{2} (1-r^{2} \rho ^{2} )(1-\left|\xi \right|^{2} {\kern 1pt} )}{(1-r^{2} \rho ^{2} )\left|1-\mathop{\xi }\limits^{\_ } \rho \zeta \right|^{2} } \rho dm\left(\zeta \right)d\rho \le   $$

\

     $$\le 2r^{2} \left\| f\right\| _{B_{\phi } }^{2} \int \limits _{T}\frac{1-\left|\mathop{\xi }\limits^{\_ } \rho \zeta \right|^{2} }{\left|1-\mathop{\xi }\limits^{\_ } \rho \zeta \right|^{2} } dm\left(\zeta \right)\int \limits _{0}^{1}\frac{\phi ^{2} (1-r^{2} \rho ^{2} {\kern 1pt} )}{1-r^{2} \rho ^{2} } \rho d\rho =  $$

\

     $$=\left\| f\right\| _{B_{\phi } }^{2} \int \limits _{1-r^{2} }^{1}\frac{\phi ^{2} (t)}{t} dt. $$

      Here  we used the identity  $$\int \limits _{T}\frac{1-\left|z\right|^{2} }{\left|1-\zeta z\right|^{2} } {\kern 1pt} {\kern 1pt} {\kern 1pt} {\kern 1pt} dm(\zeta )=1. $$

     Therefore,

     \

$$\left\| f_{r} \right\| _{BMOA}^{2} \le \left\| f\right\| _{B_{\phi } }^{2} g(1-r^{2} ).$$

\

     {\bf Corollary.}  {\it If  } $f\in {\kern 1pt} {\kern 1pt} {\kern 1pt} B{\kern 1pt} {\kern 1pt} $ {\it  then }

     \

     $$\left\| f_{r} \right\| _{BMOA} \le \left\| f\right\| _{B_{\phi } } \sqrt{\left|\log (1-r^{2} )\right|} {\kern 1pt} {\kern 1pt} {\kern 1pt} {\kern 1pt} ,{\kern 1pt} {\kern 1pt} {\kern 1pt} {\kern 1pt} {\kern 1pt} {\kern 1pt} {\kern 1pt} {\kern 1pt} 0<r<1.$$

\

     B.Korenblum [3] proved an analogous BMO estimate, applying the Garsia norm.

\

     {\bf Theorem 2.} {\it There are positive numerical constants   }
$\gamma $ {\it  and  } $M$ {\it  }{\it  such that for all}  $f\in
{\kern 1pt} {\kern 1pt} {\kern 1pt} B_{\phi } ,{\kern 1pt} {\kern
1pt} {\kern 1pt} {\kern 1pt} f(0)=0$

\

\begin{equation}
\int \limits _{T}\exp {\kern 1pt} {\kern 1pt} ({\kern 1pt} {\kern
1pt} {\kern 1pt} {\kern 1pt} \frac{\gamma {\kern 1pt} {\kern 1pt}
\left|f\left(r\zeta \right)\right|}{\left\| f\right\| _{B_{\phi }
} \sqrt{g(1-r^{2} )} } {\kern 1pt} {\kern 1pt} {\kern 1pt} {\kern
1pt} {\kern 1pt} ){\kern 1pt} {\kern 1pt} {\kern 1pt} {\kern 1pt}
dm(\varsigma ){\kern 1pt} {\kern 1pt} {\kern 1pt} {\kern 1pt} \le
M.
\end{equation}

\

     {\it Proof.}  The John-Nirenberg theorem   [1, 2]   says that there are positive constants  $c$  and  $C$   such that

     \

$$\frac{m(\zeta \in I:\left|f(\zeta )-f_{I} \right|>\lambda )}{m{\kern 1pt} (I)} \le C\exp {\kern 1pt} {\kern 1pt} {\kern 1pt} ({\kern 1pt} {\kern 1pt} {\kern 1pt} \frac{-c\lambda }{\left\| f\right\| _{BMOA} } {\kern 1pt} {\kern 1pt} {\kern 1pt} )$$

\

     for all  $f\in {\kern 1pt} {\kern 1pt} {\kern 1pt} BMOA,{\kern 1pt} {\kern 1pt} {\kern 1pt} {\kern 1pt} \lambda >0,{\kern 1pt} {\kern 1pt} {\kern 1pt} I\subseteq T.$

\

     If    $f(0)=0$  then

\

       $$f_{T} =\frac{1}{2\pi } \int \limits _{T}f(\varsigma )\left|d\zeta \right|=0 $$

\

        and

\

\begin{equation}
                                 E\left(\lambda
\right)=m\left(\zeta \in T:\left|f(r\zeta )\right|>\lambda
\right)\le C\exp {\kern 1pt} {\kern 1pt} {\kern 1pt} ({\kern 1pt}
{\kern 1pt} {\kern 1pt} \frac{-c\lambda }{\left\| f\right\|
_{BMOA} } {\kern 1pt} {\kern 1pt} {\kern 1pt} ) .
 \end{equation}

\

     Since  $E\left(\lambda \right)$   is the distilution function of $f$ ,  then for all  $p>0$ [1]

\

\begin{equation}
                                     \int \limits
_{T}\left|f\right|^{p} dm =p\int \limits _{0}^{\infty }\lambda
^{p-1} E(\lambda ){\kern 1pt} {\kern 1pt} d\lambda  .
\end{equation}

\

     If {\kern 1pt}    $0<\gamma <c$ ,  using (4) and (3), we obtain

\

$$\int \limits _{T}\exp ({\kern 1pt} {\kern 1pt} {\kern 1pt} {\kern 1pt} \frac{\gamma {\kern 1pt} {\kern 1pt} \left|f(\zeta )\right|}{\left\| f\right\| _{BMOA} } {\kern 1pt} {\kern 1pt} {\kern 1pt} {\kern 1pt} {\kern 1pt} ){\kern 1pt} {\kern 1pt} dm\left(\zeta \right)=1+ \sum \limits _{n\ge 1}\frac{1}{n!} \frac{\gamma ^{n} }{(\left\| f\right\| _{BMOA} )^{n} }  \int \limits _{T}\left|f(\zeta )\right|^{n} dm\left(\zeta \right)= $$

\

$$=1+\sum \limits _{n\ge 1}\frac{1}{n!} \frac{\gamma ^{n} }{(\left\| f\right\| _{BMOA} )^{n} }  {\kern 1pt} {\kern 1pt} {\kern 1pt} {\kern 1pt} {\kern 1pt} {\kern 1pt} {\kern 1pt} n\int \limits _{0}^{\infty }\lambda ^{n-1} E(\lambda )d\lambda  =$$

\

$$=1+\frac{\gamma }{\left\| f\right\| _{BMOA} } {\kern 1pt} {\kern 1pt} {\kern 1pt} {\kern 1pt} {\kern 1pt} {\kern 1pt} {\kern 1pt} {\kern 1pt} \sum \limits _{n\ge 1}\frac{1}{n-1!} \frac{\gamma ^{n-1} }{(\left\| f\right\| _{BMOA} )^{n-1} }  {\kern 1pt} {\kern 1pt} {\kern 1pt} {\kern 1pt} {\kern 1pt} {\kern 1pt} {\kern 1pt} \int \limits _{0}^{\infty }\lambda ^{n-1} E(\lambda )d\lambda  =$$

\

$$=1+\frac{\gamma }{\left\| f\right\| _{BMOA} } \int \limits _{0}^{\infty }E(\lambda ) \sum \limits _{n\ge 1}\frac{1}{n-1!} \frac{\gamma ^{n-1} \lambda ^{n-1} }{(\left\| f\right\| _{BMOA} )^{n-1} }  {\kern 1pt} {\kern 1pt} {\kern 1pt} {\kern 1pt} {\kern 1pt} {\kern 1pt} d\lambda =$$

\

$$=1+\frac{\gamma }{\left\| f\right\| _{BMOA} } \int \limits _{0}^{\infty }E(\lambda )\exp ({\kern 1pt} {\kern 1pt} {\kern 1pt} \frac{\gamma \lambda }{\left\| f\right\| _{BMOA} } {\kern 1pt} {\kern 1pt} {\kern 1pt} ){\kern 1pt} {\kern 1pt} {\kern 1pt} d\lambda  \le $$

\

$$\le 1+\frac{\gamma }{\left\| f\right\| _{BMOA} } {\kern 1pt} {\kern 1pt} {\kern 1pt} {\kern 1pt} C\int \limits _{0}^{\infty }\exp ({\kern 1pt} {\kern 1pt} {\kern 1pt} {\kern 1pt} \frac{-\left(c-\gamma \right)\lambda }{\left\| f\right\| _{BMOA} }
{\kern 1pt} {\kern 1pt} {\kern 1pt} {\kern 1pt} ){\kern 1pt} {\kern 1pt} d\lambda = {\kern 1pt} {\kern 1pt} {\kern 1pt} {\kern 1pt} {\kern 1pt} {\kern 1pt} 1+\frac{\gamma C}{c-\gamma } {\kern 1pt} {\kern 1pt} {\kern 1pt} {\kern 1pt} {\kern 1pt}
{\kern 1pt} {\kern 1pt} {\kern 1pt} {\kern 1pt} {\kern 1pt} {\kern 1pt} {\kern 1pt} {\kern 1pt} \mathop{=}\limits^{def} {\kern 1pt} {\kern 1pt} {\kern 1pt} {\kern 1pt} {\kern 1pt} {\kern 1pt} {\kern 1pt} {\kern 1pt} {\kern 1pt} M<\infty {\kern 1pt}
{\kern 1pt} {\kern 1pt} {\kern 1pt} .$$

\

     Putting  $f=f_{r} $  and applying (1), we obtain (2).

\

     {\bf Theorem 3.} There is a constant  $\gamma _{1} $ , such that
for every  $f\in {\kern 1pt} {\kern 1pt} {\kern 1pt} B_{\phi } $ ,
$f\left(0\right)=0$

     \

\begin{equation}
                                \mathop{\lim }\limits_{r\to
1^{-} } \sup \frac{\left|f(r\zeta )\right|}{\log \left|\log
\left(1-r\right)\right|\sqrt{g(1-r^{2} )} } \le \gamma _{1}
\left\| f\right\| _{B_{\phi } }
\end{equation}

\

     for almost all  $\zeta \in T.$

\

     {\it Proof.}  Theorem 2  implies that  $(0<r<1)$

     \

$$\int \limits _{0}^{1}{\kern 1pt} {\kern 1pt} ({\kern 1pt} {\kern 1pt} \frac{1}{\left(1-r\right)\log ^{2} {\kern 1pt} (\frac{e}{1-r} {\kern 1pt} {\kern 1pt} {\kern 1pt} )}  {\kern 1pt} {\kern 1pt} {\kern 1pt} {\kern 1pt} \int \limits _{T}\exp ({\kern 1pt} {\kern 1pt} {\kern 1pt} {\kern 1pt} \frac{\gamma \left|f(r\zeta )\right|}{\left\| f\right\| _{B_{\phi } } {\kern 1pt} {\kern 1pt} {\kern 1pt} {\kern 1pt} {\kern 1pt} {\kern 1pt} \sqrt{g(1-r^{2} )} } ){\kern 1pt} {\kern 1pt} dm(\zeta ){\kern 1pt} {\kern 1pt} {\kern 1pt} {\kern 1pt} ){\kern 1pt} {\kern 1pt} dr\le M. $$

\

     Therefore , for almost all  $\zeta $

     \

$$\int \limits _{0}^{1}({\kern 1pt} {\kern 1pt} {\kern 1pt} \frac{1}{\left(1-r\right)\log ^{2} {\kern 1pt} (\frac{e}{1-r} {\kern 1pt} {\kern 1pt} )}  {\kern 1pt} {\kern 1pt} {\kern 1pt} {\kern 1pt} {\kern 1pt} {\kern 1pt} \exp ({\kern 1pt} {\kern 1pt} {\kern 1pt} {\kern 1pt} {\kern 1pt} {\kern 1pt} {\kern 1pt} \frac{\gamma \left|f(r\zeta )\right|}{\left\| f\right\| _{B_{\phi } } {\kern 1pt} {\kern 1pt} {\kern 1pt} {\kern 1pt} {\kern 1pt} {\kern 1pt} \sqrt{g(1-r^{2} )} } ){\kern 1pt} {\kern 1pt} {\kern 1pt} dr{\kern 1pt} {\kern 1pt} {\kern 1pt} {\kern 1pt} {\kern 1pt} <\infty ,$$

\

which implies that

\

$$\mathop{\lim }\limits_{r{\kern 1pt} {\kern 1pt} {\kern 1pt} {\kern 1pt} \to {\kern 1pt} {\kern 1pt} {\kern 1pt} {\kern 1pt} {\kern 1pt} 1^{-} } \int \limits _{r}^{(r+1)/2}\frac{1}{\left(1-\rho \right)\log ^{2} {\kern 1pt} {\kern 1pt} {\kern 1pt} {\kern 1pt} {\kern 1pt} ({\kern 1pt} {\kern 1pt} \frac{e}{1-\rho } {\kern 1pt} {\kern 1pt} {\kern 1pt} {\kern 1pt} )} \exp ({\kern 1pt} {\kern 1pt} {\kern 1pt} {\kern 1pt} {\kern 1pt} {\kern 1pt} {\kern 1pt} \frac{\gamma {\kern 1pt} {\kern 1pt} \left|f(\rho \zeta )\right|}{\left\| f\right\| _{B_{\phi } } \sqrt{g(1-\rho ^{2} )} } ) {\kern 1pt} {\kern 1pt} {\kern 1pt} d\rho {\kern 1pt} {\kern 1pt} {\kern 1pt} {\kern 1pt} {\kern 1pt} {\kern 1pt} {\kern 1pt} {\kern 1pt} =0{\kern 1pt} {\kern 1pt} .$$

\

     Putting      $\mu \left(r,\zeta \right)=\min \left\{\left|f(\rho \zeta )\right|:r\le \rho \le (r+1)/2\right\}$   we get

\

$$\int \limits _{r}^{(r+1)/2}\frac{1}{\left(1-\rho \right)\log ^{2} {\kern 1pt} ({\kern 1pt} \frac{e}{1-\rho } {\kern 1pt} {\kern 1pt} )} {\kern 1pt} {\kern 1pt} {\kern 1pt} {\kern 1pt} {\kern 1pt} {\kern 1pt} {\kern 1pt} \exp ({\kern 1pt} {\kern 1pt} \frac{\gamma \left|f(\rho \zeta )\right|}{\left\| f\right\| _{B_{\phi } } {\kern 1pt} {\kern 1pt} {\kern 1pt} {\kern 1pt} {\kern 1pt} {\kern 1pt} \sqrt{g(1-\rho ^{2} )} } ) {\kern 1pt} {\kern 1pt} {\kern 1pt} d\rho {\kern 1pt} {\kern 1pt} {\kern 1pt} {\kern 1pt} {\kern 1pt} {\kern 1pt} {\kern 1pt} {\kern 1pt} \ge $$

\

$$\ge {\kern 1pt} {\kern 1pt} {\kern 1pt} {\kern 1pt} {\kern 1pt} \int \limits _{r}^{(r+1)/2}\log ^{-2} (\frac{e}{1-\rho } {\kern 1pt} {\kern 1pt} ){\kern 1pt} {\kern 1pt} {\kern 1pt} {\kern 1pt} {\kern 1pt} {\kern 1pt} {\kern 1pt} \exp ({\kern 1pt} {\kern 1pt} \frac{\gamma \left|f(\rho \zeta )\right|}{\left\| f\right\| _{B_{\phi } } {\kern 1pt} {\kern 1pt} {\kern 1pt} {\kern 1pt} {\kern 1pt} {\kern 1pt} \sqrt{g(1-\rho ^{2} )} } ) {\kern 1pt} {\kern 1pt} {\kern 1pt} d\rho {\kern 1pt} {\kern 1pt} {\kern 1pt} {\kern 1pt} {\kern 1pt} {\kern 1pt} {\kern 1pt} {\kern 1pt} \ge $$

\

$$\ge {\kern 1pt} {\kern 1pt} {\kern 1pt} {\kern 1pt} {\kern 1pt} \int \limits _{r}^{(r+1)/2}\log ^{-2} (\frac{e}{1-\rho } {\kern 1pt} {\kern 1pt} ){\kern 1pt} {\kern 1pt} {\kern 1pt} {\kern 1pt} {\kern 1pt} {\kern 1pt} {\kern 1pt} \exp ({\kern 1pt} {\kern 1pt} \frac{\gamma \left|f(\rho \zeta )\right|}{\left\| f\right\| _{B_{\phi } } {\kern 1pt} {\kern 1pt} {\kern 1pt} {\kern 1pt} {\kern 1pt} {\kern 1pt} \sqrt{g(1-\rho ^{2} )} } ) {\kern 1pt} {\kern 1pt} {\kern 1pt} d\rho {\kern 1pt} {\kern 1pt} {\kern 1pt} {\kern 1pt} {\kern 1pt} {\kern 1pt} {\kern 1pt} {\kern 1pt} \ge $$

\

$${\kern 1pt} \ge \log ^{-2} (\frac{2e}{1-r } ){\kern 1pt} {\kern 1pt} {\kern 1pt} {\kern 1pt} {\kern 1pt} \exp {\kern 1pt} {\kern 1pt} {\kern 1pt} ({\kern 1pt} {\kern 1pt} {\kern 1pt} \frac{\gamma {\kern 1pt} {\kern 1pt} {\kern 1pt} \mu (r,\zeta )}{\left\| f\right\| _{B_{\phi } } {\kern 1pt} {\kern 1pt} \sqrt{g((1-r)(3+r)/4))} } {\kern 1pt} {\kern 1pt} {\kern 1pt} {\kern 1pt} ){\kern 1pt} {\kern 1pt} {\kern 1pt} {\kern 1pt} >{\kern 1pt} {\kern 1pt} {\kern 1pt} {\kern 1pt} 0.$$

\

     We used the inequalities

     \

$$\displaystyle \log ^{-2} (\frac{e}{1-\rho } ){\kern 1pt}
{\kern 1pt} {\kern 1pt} {\kern 1pt} \ge {\kern 1pt} {\kern 1pt}
{\kern 1pt} {\kern 1pt} {\kern 1pt} \log ^{-2}
(\frac{e}{1-(r+1)/2} ){\kern 1pt} {\kern 1pt} {\kern 1pt} =\log
^{-2} ({\kern 1pt} {\kern 1pt} \frac{2e}{1-r} {\kern 1pt} {\kern
1pt} ),$$

\

     $$\displaystyle g(1-\rho ^{2} ){\kern 1pt} {\kern 1pt} {\kern 1pt} {\kern 1pt} \le {\kern 1pt} {\kern 1pt} {\kern 1pt} {\kern 1pt} {\kern 1pt} g(1-{\kern 1pt} {\kern 1pt} (r+1/2)^{2} {\kern 1pt} ){\kern 1pt} =g((1-r)(3+r)/4)).$$

\

     Then

     \

     $$\mathop{\lim }\limits_{r{\kern 1pt} {\kern 1pt} {\kern 1pt} {\kern 1pt} \to {\kern 1pt} {\kern 1pt} {\kern 1pt} {\kern 1pt} {\kern 1pt} 1^{-} } \log ^{-2} (\frac{2e}{1-\rho } ){\kern 1pt} {\kern 1pt} {\kern 1pt} {\kern 1pt} {\kern 1pt} \exp {\kern 1pt} {\kern 1pt} {\kern 1pt} ({\kern 1pt} {\kern 1pt} \frac{\gamma {\kern 1pt} {\kern 1pt} {\kern 1pt} \mu (r,\zeta )}{\left\| f\right\| _{B_{\phi } } {\kern 1pt} \sqrt{g((1-r)(3+r)/4))} } {\kern 1pt} {\kern 1pt} {\kern 1pt} ){\kern 1pt} {\kern 1pt} {\kern 1pt} {\kern 1pt} {\kern 1pt} {\kern 1pt} {\kern 1pt} =0,$$
which implies

\

$$\mathop{\lim }\limits_{r\to 1^{-} } {\kern 1pt} {\kern 1pt} {\kern 1pt} {\kern 1pt} ({\kern 1pt} {\kern 1pt} {\kern 1pt} \frac{\gamma {\kern 1pt} {\kern 1pt} {\kern 1pt} {\kern 1pt} {\kern 1pt} \gamma {\kern 1pt} {\kern 1pt} {\kern 1pt} \mu (r,\zeta )}{\left\| f\right\| _{B_{\phi } } {\kern 1pt} {\kern 1pt} \sqrt{g((1-r)(3+r)/4))} } -2\log \log \frac{2e}{1-r} {\kern 1pt} {\kern 1pt} {\kern 1pt} ){\kern 1pt} {\kern 1pt} {\kern 1pt} {\kern 1pt} {\kern 1pt} =-\infty .$$

\

     Since

     \

        $$\mathop{\lim }\limits_{r\to 1^{-} } {\kern 1pt} {\kern 1pt} {\kern 1pt} ({\kern 1pt} {\kern 1pt} \log \log \frac{1}{1-r} -\log \log \frac{2e}{1-r} {\kern 1pt} {\kern 1pt} {\kern 1pt} ){\kern 1pt} {\kern 1pt} {\kern 1pt} {\kern 1pt} =0,$$

        \

       it can be seen easily that

       \

\begin{equation}
                           \mu \left(r,\zeta \right)<\gamma _{1}
\left\| f\right\| _{B_{\phi } } {\kern 1pt} {\kern 1pt}
\sqrt{g((1-r)(3+r)/4))} {\kern 1pt} {\kern 1pt} {\kern 1pt} {\kern
1pt} {\kern 1pt} {\kern 1pt} \log \left|\log
\left(1-r\right)\right|\le
\end{equation}

\

 $\displaystyle \le
\gamma _{1} \left\| f\right\| _{B_{\phi } } .\sqrt{g((1-r^{2}
)/2))} {\kern 1pt} {\kern 1pt} {\kern 1pt} {\kern 1pt} {\kern 1pt}
{\kern 1pt} \log \left|\log \left(1-r\right)\right|$ .

\

for almost all  $\zeta ,{\kern 1pt} {\kern 1pt} {\kern 1pt} {\kern 1pt} {\kern 1pt} r$  sufficiently close to  $1$  and  $\gamma _{1} =2/\gamma .$

\

     In addition, let

\

      $\mu (r,{\kern 1pt} {\kern 1pt} {\kern 1pt} \zeta )=\left|f(r_{1} \zeta )\right|$  ,   $r\le r_{1} \le (r+1)/2$ .

\

     Then

     \

$$\left|f(r\zeta )\right|-\mu (r,\zeta )\le \int \limits _{r}^{r_{1} }\left|f'(\rho \zeta )\right|d\rho \le \left\| f\right\| _{B_{\phi } } \int \limits _{r}^{(r+1)/2}\frac{\phi {\kern 1pt} {\kern 1pt} (1-\rho ^{2} )}{1-\rho ^{2} }   {\kern 1pt} {\kern 1pt} {\kern 1pt} {\kern 1pt} d\rho =$$

\

$$=\frac{\left\| f\right\| _{B_{\phi } } }{2r} \int \limits _{(1-r)(3+r)/4}^{1-r^{2} }\frac{\phi (t)}{t} {\kern 1pt} {\kern 1pt} {\kern 1pt}  dt{\kern 1pt} {\kern 1pt} {\kern 1pt} {\kern 1pt} {\kern 1pt} \le {\kern 1pt} {\kern 1pt} {\kern 1pt} {\kern 1pt} {\kern 1pt} {\kern 1pt} \frac{\left\| f\right\| _{B_{\phi } } }{2r} \int \limits _{(1-r^{2} )/2}^{1-r^{2} }\frac{\phi (t)}{t} {\kern 1pt} {\kern 1pt} {\kern 1pt}  dt{\kern 1pt} {\kern 1pt} {\kern 1pt} {\kern 1pt} {\kern 1pt} {\kern 1pt} {\kern 1pt} {\kern 1pt} {\kern 1pt} \le {\kern 1pt} {\kern 1pt} {\kern 1pt} {\kern 1pt} $$

\

$$\le {\kern 1pt} {\kern 1pt} {\kern 1pt} {\kern 1pt} \frac{\left\| f\right\| _{B_{\phi } } }{2r} {\kern 1pt} {\kern 1pt} {\kern 1pt} {\kern 1pt} {\kern 1pt} {\kern 1pt} ({\kern 1pt} {\kern 1pt} {\kern 1pt} {\kern 1pt} {\kern 1pt} {\kern 1pt} \int \limits _{(1-r^{2} )/2}^{1-r^{2} }\frac{\phi ^{2} (t)}{t}  {\kern 1pt} {\kern 1pt} dt{\kern 1pt} {\kern 1pt} {\kern 1pt} {\kern 1pt} {\kern 1pt} {\kern 1pt} )^{1/2} {\kern 1pt} {\kern 1pt} {\kern 1pt} \le {\kern 1pt} {\kern 1pt} {\kern 1pt} {\kern 1pt} {\kern 1pt} {\kern 1pt} {\kern 1pt} \frac{\left\| f\right\| _{B_{\phi } } }{2r} {\kern 1pt} {\kern 1pt} {\kern 1pt} ({\kern 1pt} {\kern 1pt} {\kern 1pt} {\kern 1pt} {\kern 1pt} \int \limits _{(1-r^{2} )/2}^{1}\frac{\phi ^{2} (t)}{t}  {\kern 1pt} {\kern 1pt} dt{\kern 1pt} {\kern 1pt} {\kern 1pt} {\kern 1pt} {\kern 1pt} )^{1/2} ={\kern 1pt} {\kern 1pt} {\kern 1pt} {\kern 1pt} {\kern 1pt} $$

\

$$={\kern 1pt} {\kern 1pt} {\kern 1pt} {\kern 1pt} {\kern 1pt} \frac{\left\| f\right\| _{B_{\phi } } }{2r} \sqrt{g((1-r^{2} )/2))} {\kern 1pt} {\kern 1pt} {\kern 1pt} {\kern 1pt} .$$

\

     Applying (6), we obtain

\

 $\left|f(r\zeta )\right|\le $  $\left\| f\right\| _{B_{\phi } } \sqrt{g((1-r^{2} )/2))} {\kern 1pt} {\kern 1pt} {\kern 1pt} {\kern 1pt} {\kern 1pt} {\kern 1pt} {\kern 1pt} ({\kern 1pt} {\kern 1pt} \frac{1}{2r} +\gamma _{1} \log \left|\log \left(1-r\right)\right|{\kern 1pt} {\kern 1pt} {\kern 1pt} ){\kern 1pt} {\kern 1pt}
 $

 \

for almost all
 $\zeta $ and  $r$  sufficiently close to  $1$ ,

 \

which proves  (5).

\

\

     {\bf Corollary.} (Korenblum [3])  {\it If  } $f\in {\kern 1pt} {\kern 1pt} {\kern 1pt} B,{\kern 1pt} {\kern 1pt} {\kern 1pt} {\kern 1pt} {\kern 1pt} {\kern 1pt} {\kern 1pt} {\kern 1pt} f\left(0\right)=0{\kern 1pt} {\kern 1pt} {\kern 1pt} {\kern 1pt} {\kern 1pt} {\kern 1pt} {\kern 1pt} {\kern 1pt} {\kern 1pt} {\kern 1pt} $ {\it then}

\

$$\mathop{\lim }\limits_{r\to 1^{-} } \sup \frac{\left|f(r\zeta )\right|}{\sqrt{\left|\log \left(1-r\right)\right|} {\kern 1pt} {\kern 1pt} {\kern 1pt} {\kern 1pt} \log \left|\log \left(1-r\right)\right|} {\kern 1pt} {\kern 1pt} {\kern 1pt} {\kern 1pt} {\kern 1pt} \le k{\kern 1pt} {\kern 1pt} {\kern 1pt} {\kern 1pt} {\kern 1pt} \left\| f\right\| _{B} $$

\

     {\it for almost all}{\it  } $\zeta \in T$ {\it ,  }{\it where } $k$ {\it  }{\it is an absolute constant.}

\

\

\

\

\

\

\

\

\

\

\

\

\noindent
{\small Department of Mathematics\\
        Technical University\\
        25, Tsanko Dijstabanov,\\
        Plovdiv, Bulgaria\\
        e-mail: peyyyo@mail.bg}
\end{document}